%% file: neurips_2026_workshop.tex
\documentclass{article}

\usepackage[preprint]{neurips_2026}

\usepackage[utf8]{inputenc} 
\usepackage[T1]{fontenc}    
\usepackage{hyperref}       
\usepackage{url}            
\usepackage{booktabs}       
\usepackage{amsfonts}       
\usepackage{nicefrac}       
\usepackage{microtype}      
\usepackage{xcolor}         
\usepackage{graphicx}
\usepackage{amsmath}
\usepackage{subcaption}

\title{SurgeGen: A Hybrid Generative Diffusion Framework for
Storm Surge Scenario Synthesis}

\author{%
  Shunan Zheng \quad John J. Hasenbein \\
  Operations Research and Industrial Engineering\\
  University of Texas at Austin\\
  Austin, TX, USA\\
  \texttt{sn.zheng@utexas.edu} \quad
  \texttt{has@me.utexas.edu}
}

\begin{document}

\maketitle

\begin{abstract}

Predicting storm surge induced by landfalling tropical cyclones is crucial for flood mitigation and coastal risk management. Traditionally, physics-based numerical models simulate storm surge by solving the Navier--Stokes equations using numerical methods, but these simulations are computationally expensive. Generative models are promising for storm surge emulation because they can generate diverse realizations rather than producing a single deterministic prediction. However, their use for storm surge emulation remains largely unexplored. In this paper, we leverage diffusion models for storm surge surrogate modeling, combining a baseline prediction stage with conditional generation to provide a more interpretable modeling framework. We develop SurgeGen, a two-stage generative framework for generating storm surge scenarios conditioned on hypothetical storms with parameters defined in a continuous space. First, a baseline model produces a coarse estimate of the storm surge height. This estimate then conditions a diffusion model, which generates refined storm surge scenarios that better capture spatial patterns and variability. We demonstrate that our approach can generate realistic and diverse storm surge scenarios under conditions both within and outside the training distribution.
\footnote{\url{https://github.com/shunan-z/SurgeGen-framework-for-storm-surge}}

\end{abstract}

\section{Introduction}

Storm surge generated by landfalling tropical cyclones is a major source of coastal flooding and infrastructure damage, making accurate surge prediction critical for evacuation planning, coastal infrastructure design, and risk assessment. High-fidelity storm surge and coastal circulation models, including ADCIRC \cite{Luettich1992}, Delft3D \cite{Lesser2004}, and SLOSH \cite{Jelesnianski1992}, solve governing hydrodynamic equations to provide detailed simulations but are computationally expensive. Physics-based numerical models such as the Sea, Lake, and Overland Surges from Hurricanes (SLOSH) model \cite{Jelesnianski1992} simulate storm-surge responses under hypothetical hurricane scenarios, but constructing a Maximum Envelope of Water (MEOW) map requires repeatedly running SLOSH over large ensembles of hypothetical storms \cite{Glahn2009}. To support time-sensitive decision-making \cite{Austgen2025,sahin2025} before and during an approaching hurricane, models need to provide accurate scenario estimates efficiently enough for timely risk assessment  and emergency planning. Moreover, MEOW simulations are generated from discretized combinations of storm parameters, limiting efficient exploration of the continuous storm-parameter space. These computational and coverage limitations motivate data-driven surrogate models for efficient storm-surge scenario generation that generate realistic and continuous storm surge samples.

Machine learning has therefore increasingly been explored as a surrogate for storm surge and flood inundation modeling, with existing approaches focusing primarily on spatiotemporal surge prediction and peak surge estimation \cite{Qin2023}. CNNs \cite{Xie2023}, spatiotemporal recurrent models \cite{Adeli2023,Wei2024}, hierarchical representations \cite{Naeini2025}, and regional data-driven models \cite{Melsom2026} have been developed to efficiently approximate surge responses, while point-wise and neural-field approaches improve parameter efficiency and spatial generalization \cite{Pachev2023,Jiang2024} for peak surge prediction. However, these approaches typically produce a single deterministic surge map for a given set of conditions, limiting their ability to represent output variability and generate diverse scenarios for uncertainty analysis and risk assessment. Generative models provide an alternative by learning complex output distributions and enabling stochastic conditional generation. GANs \cite{Ravuri2021} and VAEs \cite{Szwarcman2024} have been used to generate spatial climate and environmental fields, while diffusion models have recently emerged as a powerful framework for generating diverse, high-dimensional physical fields \cite{Ma2024}. Diffusion-based methods have shown promising results in weather forecasting \cite{Price2025}, precipitation nowcasting \cite{Asperti2025}, and on-demand environmental scenario generation \cite{Meuer2026}, including improved ensemble reliability and the generation of diverse and extreme scenarios \cite{Asperti2025,Price2025}. Despite these advances, diffusion-based generative surrogate modeling remains largely unexplored for storm surge, where storm surge maps exhibit a distinctive zero-inflated spatial structure.

Applying generative models to storm surge emulation presents two challenges: (1) learning high-dimensional spatial surge fields from limited and expensive simulator-generated data with a storm surge inundation structure while incorporating both storm and spatial conditions, and (2) generating realistic scenarios under hypothetical conditions beyond those represented in the training data. We address these challenges with SurgeGen, a two-stage hybrid diffusion framework for peak storm surge surrogate modeling. A regression model first predicts a coarse storm surge map from storm and spatial conditions, providing an informative estimate of the surge structure. A conditional diffusion model then refines this estimate to generate detailed storm surge maps with diverse spatial patterns. SurgeGen enables efficient generation of storm-surge scenarios without repeatedly running the underlying numerical simulator and supports exploration of hypothetical storm conditions beyond the training scenarios.

Our contributions are:
\begin{itemize}
\item We develop SurgeGen, a two-stage regression--diffusion surrogate framework for conditional storm-surge scenario generation.
\item We incorporate storm characteristics and spatial conditions into the generative surrogate to capture detailed spatial surge patterns and variability.
\item We demonstrate the ability to generate realistic surge scenarios and explore conditions beyond the discrete simulator scenarios used for training.
\end{itemize}

\section{Problem formulation and methodology}

\subsection{Problem formulation}
The Sea, Lake, and Overland Surges from Hurricanes (SLOSH) model is a numerical storm-surge model that solves hydrodynamic equations on basin-specific coastal grids incorporating geographic features such as elevation and bathymetry \cite{Jelesnianski1992}. We focus on the curvilinear Galveston Bay basin and its Maximum Envelope of Water (MEOW) dataset, which contains precomputed storm surge maps for hypothetical storms defined by intensity, track direction, tide level, and forward speed \cite{noaa_meow}. Each map records the maximum water level at each grid cell, providing the simulator outputs used to train our surrogate.

We consider the problem of predicting peak surge scenarios conditioned on the tabular information of the corresponding storm. Let $x \in \mathbb{R}^{H \times W}$ denote the target peak surge scenario. Let $c_{\text{storm}} \in \mathbb{R}^{d_s}$ denote the set of global storm parameters (e.g., category and direction), and let $c_{\text{spatial}} \in \mathbb{R}^{H \times W \times d_p}$ denote bathymetrykll predictors (e.g., elevation). The storm surge map is conditioned on $c = (c_{\text{storm}}, c_{\text{spatial}})$. Our goal is to model the conditional distribution, so that it generates
\begin{equation}
\hat{x} \sim p_{\theta}(x \mid c).
\end{equation}
Direct prediction using deterministic models often fails to capture spatial uncertainty and structured errors. To address this, we come up with a two-stage framework that combines a hurdle-style hybrid model with a conditional diffusion model.

\subsection{Data preprocessing}

We use MEOW peak storm surge scenarios from the Galveston Bay basin and transform the original curvilinear SLOSH grid into a standardized $64\times64$ image representation using a nonlinear spatial transformation. The generated scenarios are transformed back to the original geographic grid for geographic interpretation and downstream analysis. Further details on the dataset and spatial transformation are provided in Appendix~\ref{sec:dataset_preprocessing}.

\subsection{Two-stage conditional model}

\begin{figure}[htbp]
    \centering
    \includegraphics[width=0.8\textwidth]{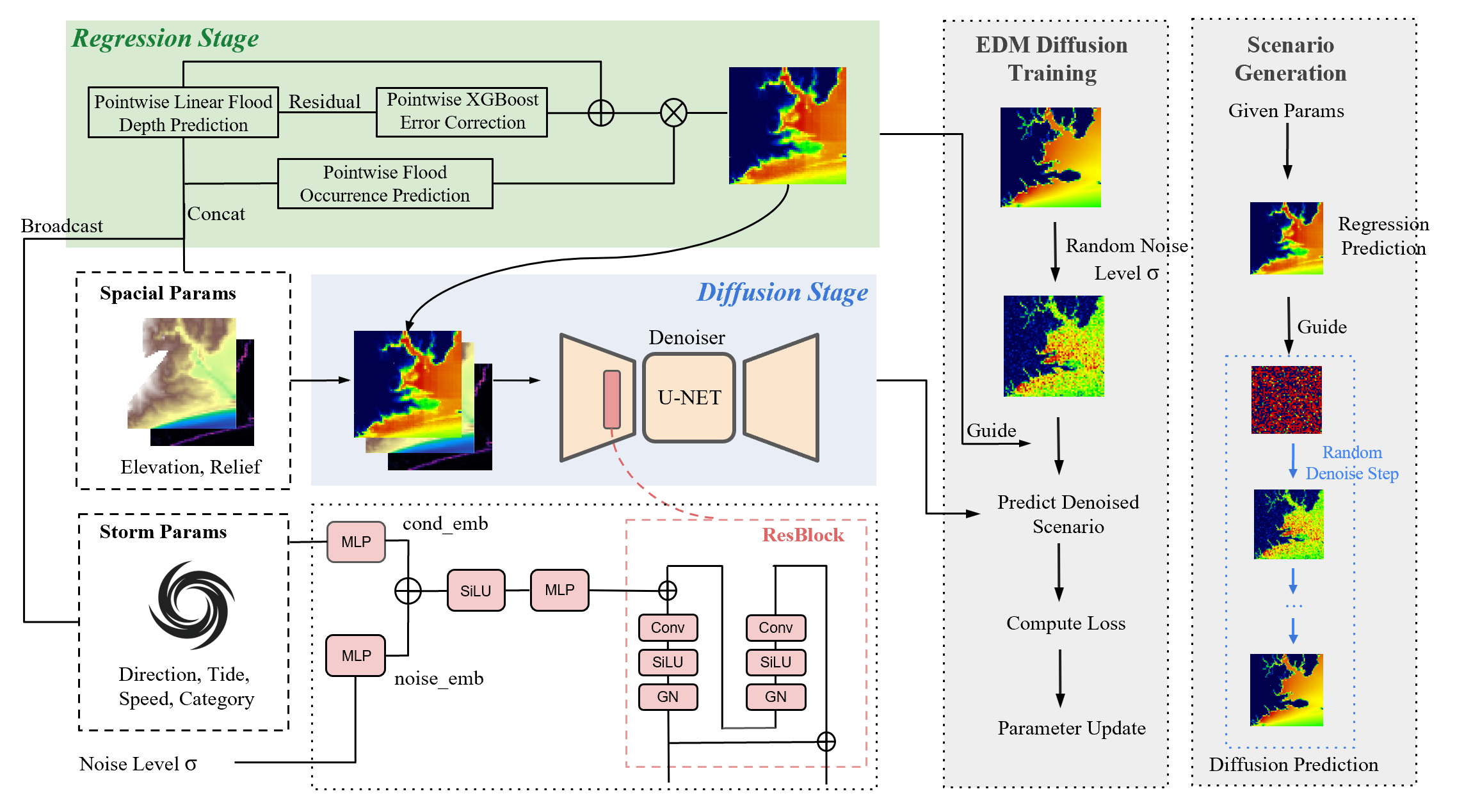}
    \caption{Our Proposed Framework and the Diffusion Training and Generation Process.}
    \label{fig:framework}
\end{figure}

We propose SurgeGen, a two-stage conditional framework for storm surge scenario generation (Figure~\ref{fig:framework}). The first stage provides a coarse point estimate of storm surge height, while the second uses SR3-style conditional diffusion \cite{Saharia2023} to generate spatially coherent refinements. Implementation details are given in Appendix~\ref{sec:model_architecture}.

\paragraph{Stage 1: Baseline prediction}

The first stage predicts a point-wise storm surge occurrence and depth from storm and spatial conditions. Because storm surge maps are strongly zero-inflated, we use a hurdle model that first predicts storm surge occurrence with logistic regression and then predicts conditional storm surge height using a combination of linear regression and gradient-boosted trees:
\begin{equation}
P(z=1\mid c)=\sigma(\beta^\top c), \qquad
\hat{x}_{\mathrm{depth}}=f_{\mathrm{lin}}(c)+f_{\mathrm{gb}}(c),
\end{equation}
where $c$ denotes the conditioning variables. The resulting coarse storm surge map assigns zero depth to locations predicted as non-inundated and the estimated depth elsewhere. This regression part provides an interpretable approximation of how surge responds to changes in storm parameters, reflecting the expected consistent variation in storm surge patterns as storm conditions change. The resulting prediction therefore provides an informative and structured condition for the diffusion stage.

\paragraph{Stage 2: Conditional diffusion}

Rather than modeling the residual between the baseline and ground truth, we use an SR3-style\cite{Saharia2023} conditional diffusion model to generate the full storm surge map conditioned on the coarse prediction. During training, Gaussian noise is added to the ground-truth map,
\begin{equation}
x_t=x_0+\sigma\epsilon,\qquad \epsilon\sim\mathcal{N}(0,I),
\end{equation}
where $\sigma$ follows the EDM noise schedule\cite{Karras2022}. The denoising network predicts the clean map $x_0$ conditioned on the coarse prediction, elevation, relief, storm parameters, and noise level. The coarse prediction provides large-scale storm surge structure, while the spatial features provide topographic context for refining local patterns.

The diffusion model uses a U-Net with multiscale convolutional blocks and skip connections. Global storm parameters are incorporated through Feature-wise Linear Modulation (FiLM) \cite{Perez2018} conditioning, where a learned storm embedding is combined with the diffusion-time embedding and injected into the residual blocks. For a target condition $c$, the first stage produces a deterministic coarse map, which guides iterative diffusion from Gaussian noise to generate the final storm surge scenario; the complete architecture and conditioning mechanism are described in Appendix~\ref{sec:model_architecture}.

\paragraph{Training objective}

The model is trained to reconstruct the clean storm surge map using a noise-weighted MSE loss with penalties for negative and small positive predictions:
\begin{equation}
\mathcal{L}=w(\sigma)
\left(
\mathcal{L}_{\mathrm{MSE}}
+\lambda_{\mathrm{small}}\mathcal{L}_{\mathrm{small}}
+\lambda_{\mathrm{neg}}\mathcal{L}_{\mathrm{neg}}
\right).
\end{equation}
The additional penalties encourage non-negative outputs and sharper separation between inundated and non-inundated regions.
\section{Experiments}

We evaluate SurgeGen on held-out MEOW scenarios and compare generated storm surge maps with the corresponding SLOSH outputs. We assess magnitude and spatial accuracy using volume error, peak error, IoU, and F1 score, and evaluate distributional similarity at infrastructure locations using Earth Mover's Distance (EMD) and Pearson correlation. Metric definitions are provided in Appendix~\ref{sec:metric_definitions}.

\subsection{Generation Quality}
\begin{figure*}[htbp]
    \centering
    \includegraphics[width=0.8\textwidth]{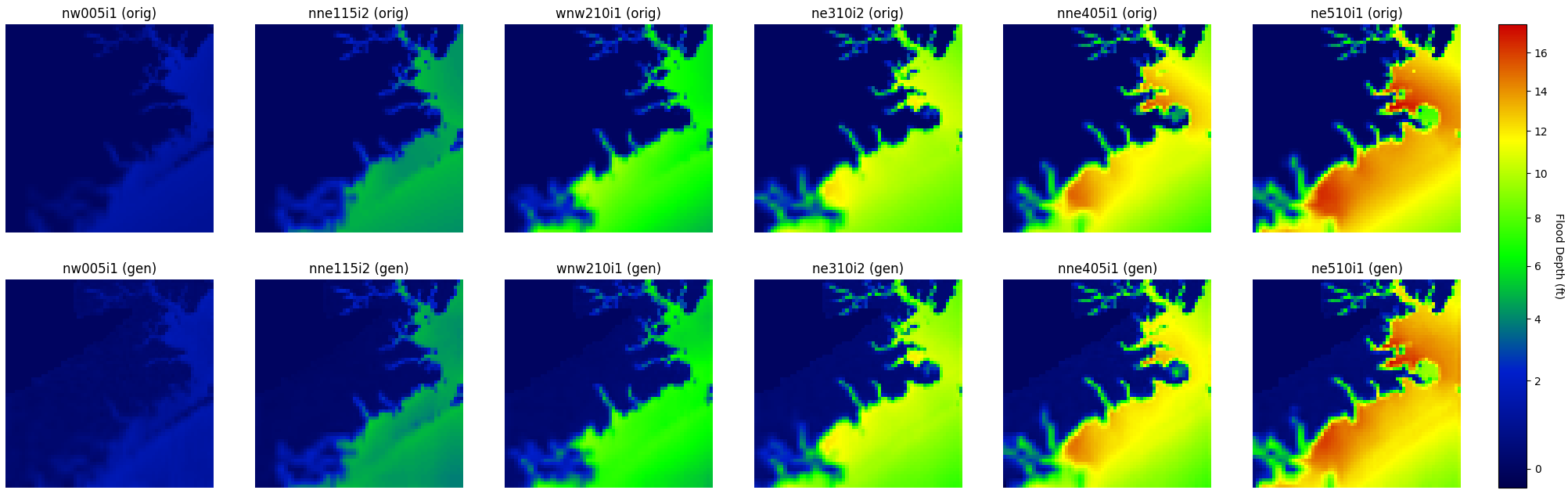}
    \caption{Comparison of test SLOSH storm surge maps with
corresponding SurgeGen samples.}
    \label{fig:flood_map}
\end{figure*}

Figure~\ref{fig:flood_map} shows representative generated storm surge scenarios. The generated scenarios visually resembles the original storm surge scenario in the first row, and inundation severely is consistently increasing as the intensity of the storm increases. Quantitatively (See Appendix~\ref{sec:metric_definitions} for metrics definitions), the model achieves a mean volume error of 0.057 ft and mean peak error of 0.484 ft. IoU and F1 scores are strongest at intermediate storm surge height thresholds, with lower scores at extreme thresholds. In practice, mitigation decisions are most important for moderate storm surge heights, where small differences in predicted severity can meaningfully affect the appropriate response.

For the downstream power-grid application, we compare generated storm surge height distribution at power-station locations with the original scenarios and the NORTA-based method \cite{Shukla2025}. SurgeGen has a higher marginal Earth Mover's Distance (EMD) than NORTA, but a lower correlation error, indicating better preservation of spatial dependence; additional distributional results are reported in Appendix~\ref{sec:additional_results}.

\subsection{Out-of-Distribution Generation}
\begin{figure*}[t]
    \centering

    \begin{minipage}[t]{0.48\textwidth}
        \centering
        \includegraphics[width=\textwidth]{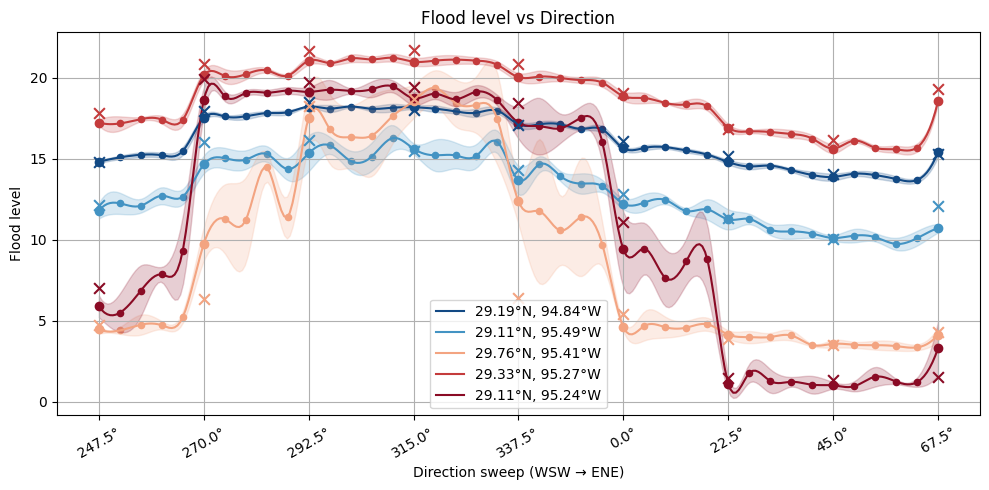}
        \captionof{figure}{Generated storm surge responses under varying directions while fixing other storm conditions.}
        \label{fig:direction}
    \end{minipage}
    \hfill
    \raisebox{3cm}{%
    \begin{minipage}[t]{0.48\textwidth}
        \centering
        \captionof{table}{Ablation study on different frameworks.}
        \label{tab:results}

        \small
        \setlength{\tabcolsep}{4pt}
        \renewcommand{\arraystretch}{0.9}

        \begin{tabular}{llcc}
            \toprule
            \textbf{Category} & \textbf{Model} & \textbf{ID} & \textbf{OOD} \\
            \midrule
            Single-stage
            & Diffusion            & 2.574 & 7.469 \\[0.5pt]
            \midrule
            Two-stage
            & XGBoost (Stage 1)    & 0.796 & 20.00 \\
            & Diffusion (with Aug.) & 1.495 & 8.039 \\
            & Residual Diffusion   & 0.854 & 7.895 \\
            & SurgeGen             & \textbf{0.751} & \textbf{5.050} \\
            \bottomrule
        \end{tabular}
    \end{minipage}%
    }

\end{figure*}

We evaluate generalization to unseen storm conditions by training on Categories 0--4 and testing on Category 5. SurgeGen generates Category 5 storm surge maps with a mean RMSE of 2.116~\(\mathrm{ft}\). We further vary the storm direction from WSW to ENE while holding other conditions fixed (Category 5 intensity, forward speed of 10 mph, and low tide level), including both directions represented in the training samples and intermediate directions. As shown in Figure~\ref{fig:direction}, the generated responses (dots) vary smoothly with direction and remain consistent with the observed scenarios (crosses), demonstrating the ability to explore conditions beyond the discrete training configurations.

\subsection{Ablation Study}
Table~\ref{tab:results} compares SurgeGen with four ablations (See Appendix~\ref{sec:additional_results} for specifications of these ablations.) that isolate the contributions of the two-stage design, the first-stage guidance, and full-map diffusion refinement. All MSE values are reported in $\mathrm{ft}^2$. In our experiment, in-distribution (ID) uses a random train/test split across Categories 0--5, while out-of-distribution (OOD) trains on Categories 0--4 and holds out Category 5 entirely for testing. 

The single-stage diffusion baseline, which removes the first-stage guidance and directly generates the storm surge map from the conditioning variables, achieves an MSE of 2.574 for ID scenarios and degrades substantially to 7.469 on OOD scenarios. Replacing the baseline prediction stage with XGBoost yields competitive ID performance (0.796) but performs poorly under OOD conditions (20.00), highlighting the importance of the baseline guidance for extrapolation. Adding stochastic augmentation increases the ID MSE to 1.495 and OOD MSE to 8.039, likely hurting because perturbing the conditioning can blur the learned storm-to-surge relationship instead of improving robustness. Finally, replacing full-map diffusion with residual diffusion results in MSEs of 0.854 and 7.895 for ID and OOD scenarios, respectively, indicating that directly refining the full storm surge map provides a more effective generative formulation than modeling the residual alone. SurgeGen achieves the lowest MSE in both settings, with 0.751 for ID and 5.050 for OOD generation. The ablation study clearly shows the contribution of each components we proposed.

\section{Limitations and Future Work}
Several limitations of this work suggest directions for future research. Our framework does not explicitly incorporate physical constraints, as doing so is challenging when peak surge is defined as the maximum over simulated storm surge maps. Richer geographic and physical information, such as coastline geometry and land cover, could further improve physical fidelity and robustness. Although our model is capable of capturing storm surge uncertainty, the current dataset does not capture the natural variability present in real storm events. Expanding the dataset by introducing purposefully designed, physically meaningful perturbations into the numerical simulations is therefore an promising direction to resolve data limitation and represent natural variability, although generating additional scenarios would be computationally expensive.
\section{Conclusion}

We propose SurgeGen, a two-stage generative surrogate for storm-surge scenario generation that combines regression-based baseline coarse prediction with conditional diffusion refinement. The framework generates spatially coherent and diverse peak storm surge scenarios while avoiding repeated numerical simulation, and experiments demonstrate realistic generation under both in-distribution and out-of-distribution conditions, including continuously varying storm directions. These results highlight the potential of diffusion-based surrogate models for efficient storm-surge scenario exploration and downstream coastal risk applications such as evacuation planning, infrastructure assessment, and flood mitigation.
\clearpage
\bibliographystyle{plain}
\bibliography{literature}

\clearpage
\appendix

\section{Dataset and preprocessing}
\label{sec:dataset_preprocessing}
\subsection{MEOW dataset details}

We use the Maximum Envelope of Water (MEOW) dataset generated by the Sea, Lake, and Overland Surges from Hurricanes (SLOSH) model for the Galveston Bay basin. The dataset contains 324 storm surge scenarios generated from combinations of four storm descriptors: hurricane intensity (Category 0--5), storm-track direction, initial tide level, and forward speed. The direction is represented by discrete approach directions spanning WSW to ENE, while tide level is represented by the low-tide ($i1$) and high-tide ($i2$) conditions. Forward speed takes values of 5, 10, and 15 mph. For each storm configuration, the SLOSH simulation provides the maximum water level at each grid cell in the basin.

Each scenario can be identified by a compact label encoding its storm conditions. For example, ``ne510i1'' denotes a Category 5 storm approaching from the northeast, with a forward speed of 10 mph and low tide. The original SLOSH basin is defined on a curvilinear grid whose spatial density varies across the basin, with finer representation near the coastline. The target storm surge maps used by SurgeGen are derived from this original grid after the spatial transformation described below.

For spatial conditioning, we use topo-bathymetric information for the Texas and Gulf of Mexico region. The transformed representation uses elevation-derived spatial features aligned to the $64\times64$ image grid. The precise construction of the spatial feature channels is kept consistent with the preprocessing used in the model implementation.

\subsection{Data preprocessing}

\begin{figure}[htbp]
    \centering
    \begin{subfigure}[t]{0.48\textwidth}
        \centering
        \includegraphics[width=\textwidth]{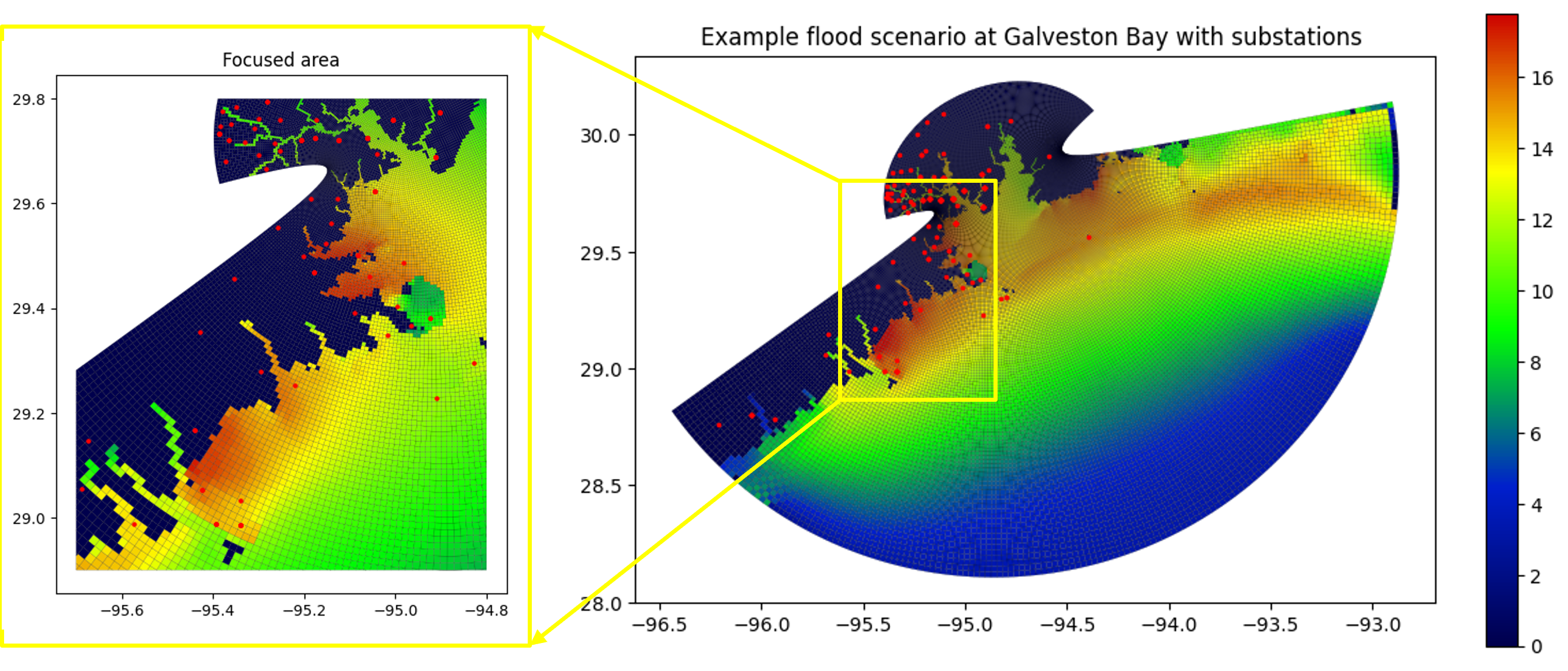}
        \caption{An example storm surge scenario at Galveston Bay. Red dots indicate substations, and the left panel shows the zoomed-in region containing substations inundated in at least one scenario.}
        \label{fig:crop}
    \end{subfigure}
    \hfill
    \begin{subfigure}[t]{0.48\textwidth}
        \centering
        \includegraphics[width=\textwidth]{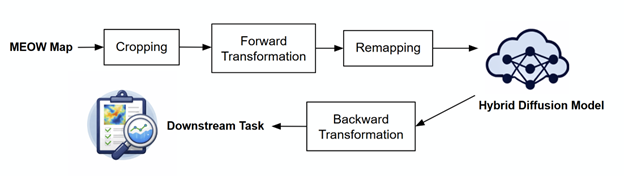}
        \caption{Data processing workflow.}
        \label{fig:data_processing}
    \end{subfigure}
    \caption{Data Processing}
\end{figure}

We first crop the Galveston Bay region using a geographic bounding box spanning $95.7^\circ$W--$94.8^\circ$W and $28.9^\circ$N--$29.8^\circ$, see Figure \ref{fig:crop}. This region contains the populated areas and infrastructure of interest and covers the area that storm surge occurs in any of the selected scenarios.

The original SLOSH grid is curvilinear and has nonuniform spatial density. In particular, coastal regions contain more and smaller cells than offshore regions. To obtain a representation compatible with convolutional diffusion models, we rescale the spatial grid so that grid cells have approximately uniform visual size and map the transformed domain to a regular $64\times64$ image grid. A nonlinear forward transformation is applied to the longitude--latitude coordinates $(x,y)$ within the cropped domain:
\begin{equation}
\tilde{x} = x\left(1 + \frac{2(x-x_1)^3}{x_1}\right),
\qquad
\tilde{y} = y\left(1 + \frac{2(y-y_1)^3}{y_1}\right)\frac{x_2-x_1}{y_2-y_1},
\end{equation}
where $(x_1,x_2)$ and $(y_1,y_2)$ define the longitude and latitude bounds of the cropped domain. The resulting raster has a nominal spatial resolution of approximately $0.014^\circ$ per pixel, corresponding to roughly 1.5 km in the study region.

For inverse reconstruction after generation, we use pull-based resampling. Rather than explicitly deriving an analytic inverse of the nonlinear transformation, we construct the desired geographic output grid, apply the known forward transformation to its coordinates, and evaluate the generated $64\times64$ image at the transformed locations using interpolation. This produces a reconstructed storm surge map aligned with the original geographic domain.

To reduce overfitting in the limited-data setting, stochastic Gaussian perturbations are applied during diffusion training. With probability $0.5$, the unperturbed input is used; otherwise, small Gaussian perturbations are introduced into the continuous spatial inputs and the continuous representations of otherwise discretized conditioning variables. The perturbation magnitude is part of the implementation and is reported separately in the training configuration. See the data processing workflow in Figure \ref{fig:data_processing}.

\section{Complete model architecture + conditioning}
\label{sec:model_architecture}
\subsection{First-stage regression}

The first stage produces a point estimate of the storm surge map from storm and spatial features and is used only as coarse guidance for the diffusion model. Because the target storm surge maps contain a large fraction of zero-valued cells, the predictor is formulated as a hurdle-style model that separates storm surge occurrence from positive storm surge height.

For each grid cell, the conditioning vector is constructed from storm descriptors and spatial features. The implementation uses elevation statistics together with storm category, forward speed, tide, and a sine/cosine representation of storm direction; geographic coordinates are also used by the pointwise first-stage model. Features are standardized or transformed according to the implementation used in training.

The first component is a logistic regression model for storm surge occurrence. Let $z\in\{0,1\}$ indicate whether a grid cell is inundated and let $c$ denote the local conditioning vector. Then
\begin{equation}
P(z=1\mid c)=\sigma(\beta^\top c),
\end{equation}
where $\sigma(\cdot)$ is the logistic function.

Conditional on storm surge occurrence, a linear regression model first estimates the storm surge height. A gradient-boosting model is then used to learn the remaining nonlinear structure and correct the linear prediction. Denoting the linear estimate by $f_{\mathrm{lin}}(c)$ and the gradient-boosting correction by $f_{\mathrm{gb}}(c)$, the conditional depth estimate is
\begin{equation}
\hat{x}_{\mathrm{depth}}=f_{\mathrm{lin}}(c)+f_{\mathrm{gb}}(c).
\end{equation}
The final baseline map assigns zero to locations classified as non-inundated and the predicted positive depth to inundated locations:
\begin{equation}
\hat{x}_{\mathrm{reg}}=
\begin{cases}
\hat{x}_{\mathrm{depth}}, & z=1,\\
0, & z=0.
\end{cases}
\end{equation}

This decomposition allows the first stage to represent the large-scale storm surge extent and magnitude while leaving spatially structured residual variation and fine-scale details to the diffusion model. The exact fitted hyperparameters of the linear, logistic, and gradient-boosting components are intentionally left for completion in the final experimental configuration.

\subsection{Diffusion model}

The second stage generates the full storm surge map rather than modeling a residual map. The clean target is \(\mathbf{x}_0\in\mathbb{R}^{1\times64\times64}\). During training, Gaussian noise with level $\sigma$ is added according to
\begin{equation}
\mathbf{x}_t=\mathbf{x}_0+\sigma\boldsymbol{\epsilon},
\qquad
\boldsymbol{\epsilon}\sim\mathcal{N}(\mathbf{0},I).
\end{equation}
The denoiser is trained to recover $\mathbf{x}_0$ from the noisy map together with the coarse regression guidance and spatial conditioning.

The denoising backbone is a three-level encoder--decoder U-Net with base width 32. The channel and spatial dimensions follow
\begin{equation}
(5,64,64)\rightarrow(32,64,64)\rightarrow(64,32,32)\rightarrow(128,16,16)
\rightarrow(256,8,8)
\end{equation}
followed by the corresponding decoder path
\begin{equation}
(256,8,8)\rightarrow(128,16,16)\rightarrow(64,32,32)
\rightarrow(32,64,64)\rightarrow(1,64,64).
\end{equation}

Each downsampling block contains two residual blocks followed by a stride-2 convolution. Each upsampling block uses a transposed convolution, concatenation with the encoder skip feature, and two residual blocks. The bottleneck contains two residual blocks at 256 channels. Residual blocks use Group Normalization with 8 groups, SiLU activations, and $3\times3$ convolutions, with a $1\times1$ projection in the skip path when the input and output channel dimensions differ.

At the tensor level, the diffusion network uses the noisy storm surge map, the first-stage regression map, elevation-derived spatial features, and the valid-cell mask as its spatial input. The exact channel-by-channel raster construction is retained from the implementation and is summarized here at the technical level rather than repeating the higher-level verbal description in the main paper.

\subsection{Conditioning and data augmentation}

The global storm condition is represented as
\begin{equation}
\mathbf{c}_{\mathrm{global}} =
[\mathrm{Category},\mathrm{Speed},\mathrm{Tide},
\mathrm{Direction}~sin,\mathrm{Direction}~cos]
\in \mathbb{R}^{5}.
\end{equation}
The storm condition is mapped to a 128-dimensional embedding using a two-layer multilayer perceptron with dimensions
\begin{equation}
5\rightarrow512\rightarrow128.
\end{equation}

The diffusion noise level is represented by
\begin{equation}
c_{\mathrm{noise}}=0.25\log\sigma.
\end{equation}
This scalar is encoded using a sinusoidal embedding of dimension 128 and passed through a second MLP with architecture
\begin{equation}
128\rightarrow512\rightarrow128,
\end{equation}
producing the time embedding $\mathbf{e}_t$. The global storm variables produce $\mathbf{e}_c$ through the corresponding conditioning MLP. The embeddings are fused by elementwise summation,
\begin{equation}
\mathbf{e}=\mathbf{e}_t+\mathbf{e}_c.
\end{equation}

For each residual block, the fused embedding is mapped by a learned linear projection to the channel dimension. If $\mathbf{h}$ denotes the feature map after the first convolution, the conditioning vector is transformed as
\begin{equation}
\mathbf{b}=W_e\,\mathrm{SiLU}(\mathbf{e})+\mathbf{a},
\qquad
\mathbf{b}\in\mathbb{R}^{C_{\mathrm{out}}},
\end{equation}
and broadcast across spatial locations:
\begin{equation}
h'_{c,i,j}=h_{c,i,j}+b_c.
\end{equation}
The block then applies the second GroupNorm, SiLU activation, and convolution before adding the residual skip connection. This is referred to as FiLM-style conditioning because the external variables modulate feature channels throughout the U-Net; technically, the present implementation uses additive channel-wise modulation rather than the full affine scale-and-shift FiLM transformation.

\input{config2.tex}

\section{Metric definitions}
\label{sec:metric_definitions}

We use complementary metrics to evaluate magnitude accuracy, spatial agreement, and distributional similarity.

\paragraph{Volume error.}
For a generated storm surge map $\hat{\mathbf{x}}$ and the corresponding reference map $\mathbf{x}$, the total storm surge volume is obtained by summing storm surge height over the valid spatial cells, with the corresponding spatial-area weighting when converting the raster values to volume. We report the absolute relative discrepancy between generated and reference volume:
\begin{equation}
E_{\mathrm{vol}}
=\frac{|V(\hat{\mathbf{x}})-V(\mathbf{x})|}{V(\mathbf{x})},
\end{equation}
where $V(\cdot)$ denotes the total storm surge volume. 

\paragraph{Peak error.}
Peak error measures the discrepancy between the maximum storm surge heights of the generated and reference scenarios:
\begin{equation}
E_{\mathrm{peak}}=|\max_i\hat{x}_i-\max_i x_i|.
\end{equation}
We report the mean peak error across evaluated scenarios.

\paragraph{Intersection over Union and F1.}
For a chosen storm-surge-height threshold, each map is converted to a binary inundated-region mask. Let $P$ and $G$ denote the predicted and ground-truth inundated regions. Then
\begin{equation}
\mathrm{IoU}=\frac{|P\cap G|}{|P\cup G|},
\qquad
F_1=\frac{2PR}{P+R},
\end{equation}
where $P$ and $R$ in the second expression denote precision and recall, respectively. The experiments use relative thresholds defined with respect to the maximum ground-truth storm surge height.

\paragraph{Mean squared error.}
For a generated map and its reference map, the map-level MSE is
\begin{equation}
\mathrm{MSE}=\frac{1}{|\Omega|}\sum_{i\in\Omega}(\hat{x}_i-x_i)^2,
\end{equation}
where $\Omega$ denotes the set of valid grid cells. For multiple generated samples, the reported value is averaged over the evaluated samples and scenarios.

\paragraph{Earth Mover's Distance.}
EMD is used to compare the empirical distributions of storm surge levels at infrastructure locations. Let $F$ and $G$ denote the empirical distributions of generated and reference storm surge levels; the reported EMD is the minimum transport cost required to transform one empirical distribution into the other under the chosen one-dimensional ground metric. 

\paragraph{Pearson correlation and correlation error.}
For infrastructure locations or scenarios, Pearson correlation measures linear dependence between generated and reference storm surge levels:
\begin{equation}
\rho_{XY}=\frac{\operatorname{Cov}(X,Y)}{\sigma_X\sigma_Y}.
\end{equation}
The downstream table reports correlation error rather than correlation itself. 

\section{Additional results}
\label{sec:additional_results}

The following results are not reproduced in the main body of the current NeurIPS version due to paper limit.

\subsection{Threshold-based spatial agreement}

\begin{figure}[ht]
    \centering
    \includegraphics[width=0.4\columnwidth]{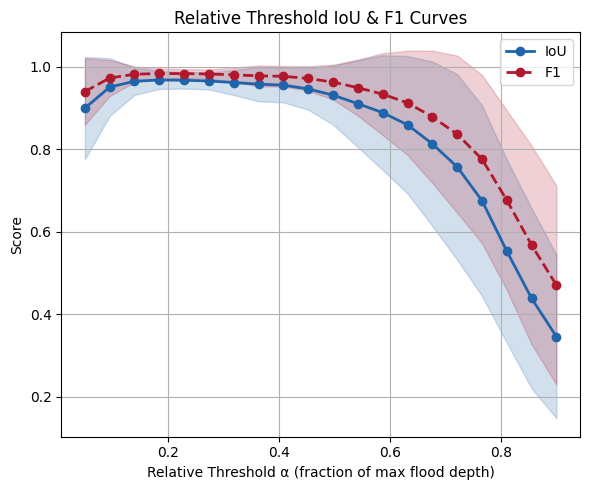}
    \caption{IoU and F1 scores across relative storm-surge-height thresholds. The curves peak at intermediate thresholds, while the uncertainty intervals reflect variation across five generated samples per scenario.}
    \label{fig:additional_iou_f1}
\end{figure}

We evaluated IoU and F1 across a range of relative storm-surge-depth thresholds. The resulting curves exhibit a bell-shaped pattern (Figure \ref{fig:additional_iou_f1}), with strongest spatial agreement at intermediate thresholds. At low thresholds, small nonzero prediction differences can substantially enlarge the predicted inundated region and reduce precision. At high thresholds, the inundated region becomes smaller and more sensitive to local spatial misalignment. Each scenario was evaluated using five generated samples, and the uncertainty bands were computed across these samples.

\begin{figure*}[ht]
    \centering
    \includegraphics[width=0.8\textwidth]{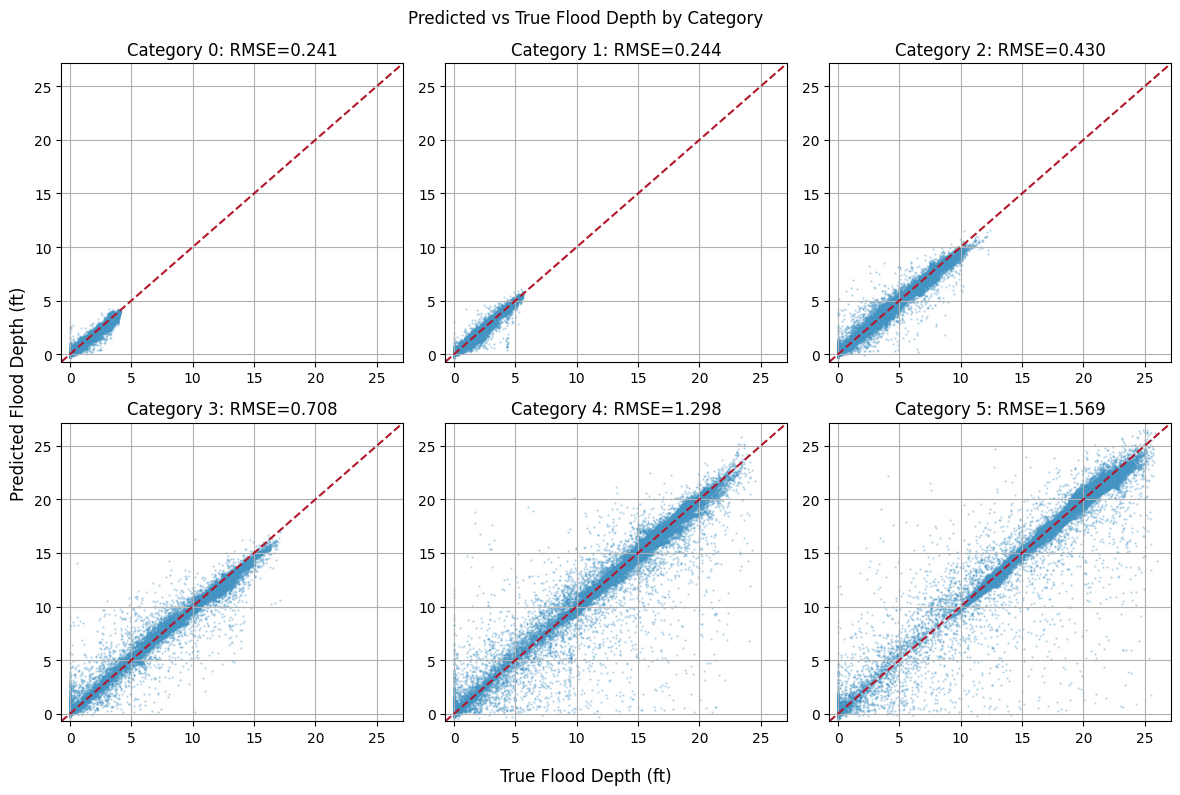}
    \caption{Predicted storm surge height versus true storm surge height for in-distribution held-out set.}
    \label{fig:scatter error id}
\end{figure*}

See Figure \ref{fig:scatter error id} for the error between predicted and true storm surge height by direction.

\subsection{Out-of-distribution depth agreement}
\label{appendix: ood}
\begin{figure*}[ht]
    \centering
    \includegraphics[width=0.8\textwidth]{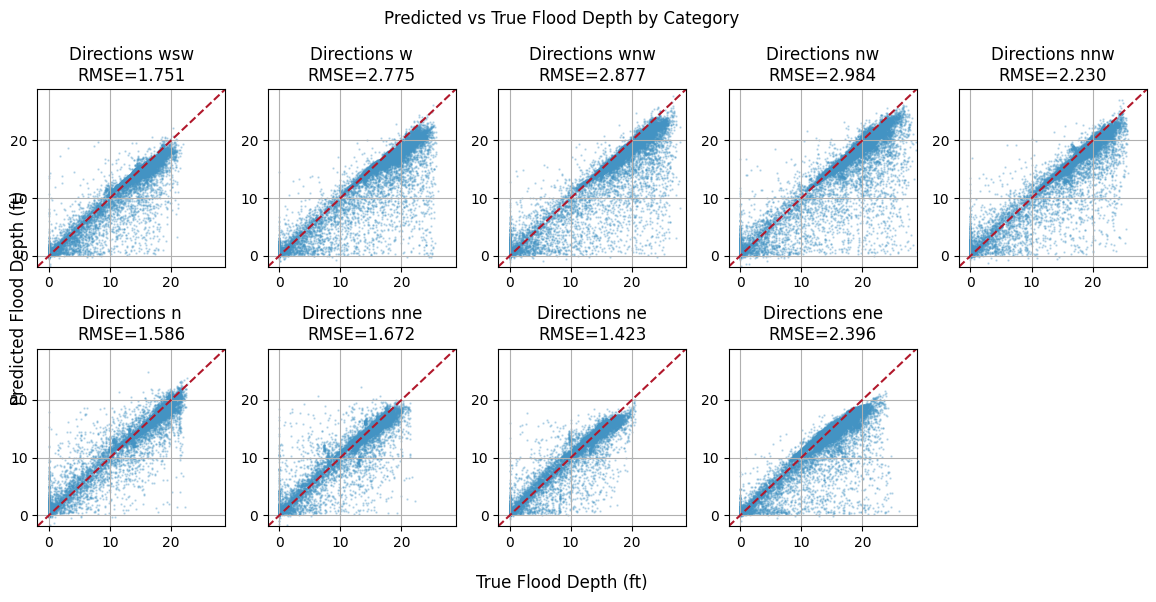}
    \caption{Predicted storm surge height versus true storm surge height for Category 5 extrapolation.}
    \label{fig:additional_scatter_ood}
\end{figure*}

\begin{figure*}[htbp]
    \centering
    \includegraphics[width=0.7\textwidth]{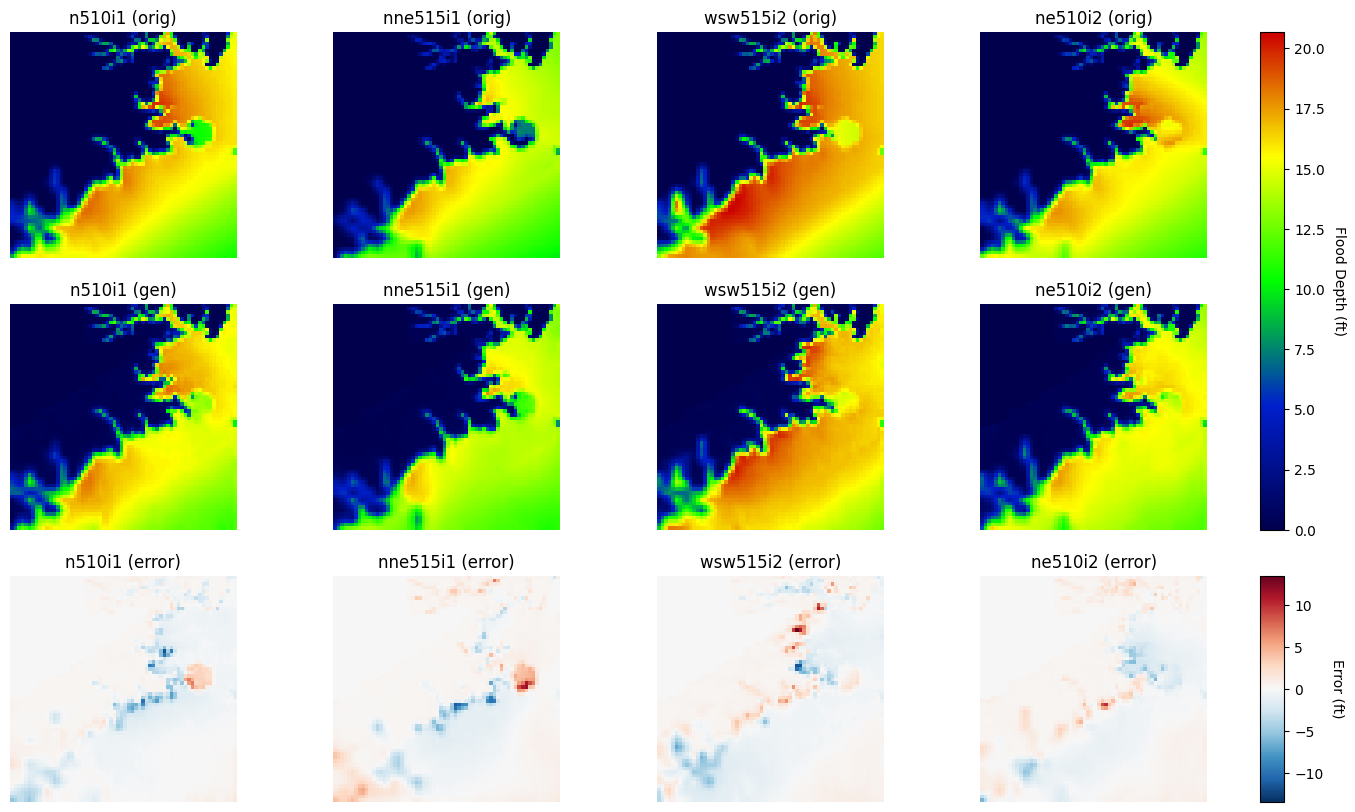}
    \caption{Example out-of-distribution generated samples for category 5 hurricane. The first row shows the original maps, the second row shows generated samples and the third row shows the error . For some scenarios, the model generate synthetic samples that match the ground truth storm surge map.}
    \label{fig:ood flood map}
\end{figure*}

For category extrapolation, the model was trained on Categories 0--4 and evaluated on Category 5 scenarios. The generated storm surge maps were compared with the corresponding SLOSH maps at the cell level. The scatter plot (Figure \ref{fig:additional_scatter_ood}) below provides an additional view of the agreement between predicted and true storm surge heights under this out-of-distribution setting. Figure \ref{fig:ood flood map} gives some examples of OOD generation result.

\subsection{Geographic visualization of generated scenarios}

\begin{figure}[ht]
    \centering
    \includegraphics[width=0.5\columnwidth]{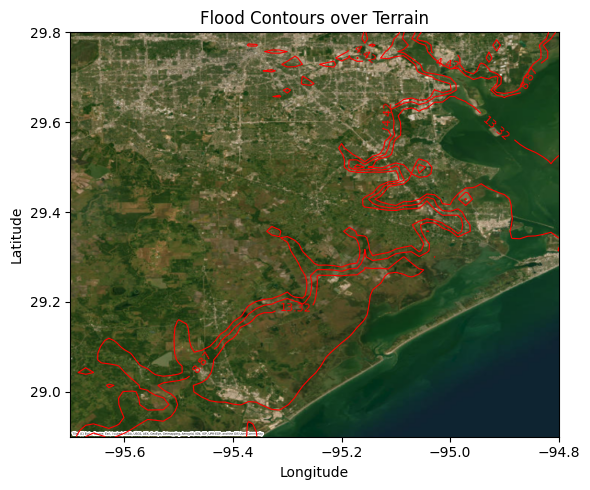}
    \caption{Geographic overlay and contour visualization of a generated storm surge scenario after mapping the model output back to the original geographic domain.}
    \label{fig:additional_contour}
\end{figure}

The earlier study also mapped generated scenarios back to the geographic domain and compared them with a satellite basemap. The corresponding contour visualization (Figure \ref{fig:additional_contour}) illustrates that the generated storm surge field preserves the expected geographic concentration of inundation near the coast and the spatial decay of storm surge height farther inland.

\subsection{Downstream power-grid comparison with NORTA}

We evaluated generated storm surge levels at power-station locations and compared the proposed diffusion model with a NORTA-based scenario-generation method (Table \ref{tab:additional_norta}). NORTA achieved lower marginal-distribution error according to EMD, whereas the diffusion model achieved lower correlation error, indicating stronger agreement in the spatial dependence structure. The full descriptive statistics from the earlier experiment are reproduced below. 

\begin{table}[ht]
\centering
\caption{Comparison of EMD and correlation error statistics for NORTA and the diffusion model in the power-grid downstream task.}
\label{tab:additional_norta}
\small
\begin{tabular}{lcccc}
\toprule
& \multicolumn{2}{c}{\textbf{NORTA}} & \multicolumn{2}{c}{\textbf{Diffusion}} \\
\textbf{Statistic} & \textbf{EMD} & \textbf{Corr. Err.} & \textbf{EMD} & \textbf{Corr. Err.} \\
\midrule
Mean & 0.2772 & 0.1139 & 0.7135 & 0.0686 \\
Std & 0.2003 & 0.1190 & 0.3974 & 0.0812 \\
Min & 0.0003 & 0.0000 & 0.0910 & 0.0000 \\
25\% & 0.1360 & 0.0179 & 0.3927 & 0.0157 \\
50\% & 0.1989 & 0.0735 & 0.6286 & 0.0377 \\
75\% & 0.4122 & 0.1805 & 0.9953 & 0.0883 \\
Max & 0.8510 & 0.7331 & 1.7969 & 0.5986 \\
\bottomrule
\end{tabular}\end{table}

\subsection{Ablation Settings}

\begin{itemize}
    \item \textbf{Single-stage Diffusion:} Removes the first-stage baseline model and directly generates the storm surge map conditioned on the storm and spatial variables.
    
    \item \textbf{XGBoost-guided:} Replaces the baseline model with XGBoost while keeping the conditional diffusion refinement unchanged.
    
    \item \textbf{Diffusion with Augmentation:} Adds the stochastic Gaussian perturbation during diffusion training to evaluate its effect on generalization. With probability $0.5$, the original continuous inputs are retained; otherwise, small Gaussian perturbations are applied to continuous spatial features and to conditioning variables representing quantities that are discretized in the original MEOW scenarios.
    
    \item \textbf{Residual Diffusion:} Uses the two-stage framework but trains the diffusion model to predict the residual between the baseline model output and the ground-truth storm surge map instead of generating the full storm surge map.
    
    \item \textbf{SurgeGen:} Uses the full proposed framework, combining baseline model guidance, stochastic augmentation, and full-map SR3-style diffusion refinement.
\end{itemize}

For all variants, we evaluate both in-distribution (ID) and out-of-distribution (OOD) performance. The ID setting uses a random split across Categories 0--4, while the OOD setting trains on Categories 0--4 and holds out Category 5 for testing. All other experimental procedures are kept consistent across variants.

\clearpage

\end{document}

%% file: config2.tex
%
%

\section{Complete Training Hyperparameters}
\label{app:hyperparameters}

Table~\ref{tab:hyperparameters} lists the full configuration used for the
reported model. The final epoch-100 checkpoint is used without checkpoint selection.

\begin{table}[htbp]
\centering
\caption{Complete training and sampling hyperparameters. Flood depths are in
metres, and $\sigma_{x_0}$ denotes the standard deviation of the training targets
over valid cells, so that the EDM noise scales
are invariant to the units of depth.}
\label{tab:hyperparameters}
\begin{tabular}{ll}
\toprule
\multicolumn{2}{l}{\textit{Data}} \\
Grid resolution            & $64 \times 64$ \\
Scenarios (train / test)   & 259 / 65, random scenario-level $80/20$ split, seed 42 \\
Spatial conditioning       & Regression baseline, mean elevation, elevation variance \\
Global conditioning        & Category, forward speed, tide, $\sin\theta$(direction), $\cos\theta$(direction) \\
\midrule
\multicolumn{2}{l}{\textit{Denoiser}} \\
Architecture               & U-Net, 3 resolutions ($64 \to 8$), 2 residual blocks per stage \\
Channel widths             & 32, 64, 128, 256 \\
Normalisation, activation  & GroupNorm (8 groups), SiLU; no attention \\
Conditioning               & Sinusoidal noise embedding (dim 128) plus global MLP, summed \\
Parameters                 & 8.86\,M \\
\midrule
\multicolumn{2}{l}{\textit{EDM noise}} \\
$\sigma_{\mathrm{data}}$   & $\sigma_{x_0}$ \\
$[\sigma_{\min}, \sigma_{\max}]$ & $[0.004\,\sigma_{x_0},\ 160\,\sigma_{x_0}]$ \\
$P_{\mathrm{mean}}$, $P_{\mathrm{std}}$ & $-1.2 + \log(\sigma_{x_0}/0.5)$,\ \ $1.2$ \\
Loss weighting             & $\lambda(\sigma) = (\sigma^2 + \sigma_{\mathrm{data}}^2)/(\sigma\sigma_{\mathrm{data}})^2$ \\
\midrule
\multicolumn{2}{l}{\textit{Objective}} \\
Reconstruction             & Masked MSE over valid cells \\
Negative-depth penalty     & $\lambda_{\mathrm{neg}} = 10$, dead band $\epsilon = 0.1$~ft \\
Shallow-depth penalty      & $\lambda_{\mathrm{small}} = 0.1$, threshold $\tau = 0.5$~ft \\
\midrule
\multicolumn{2}{l}{\textit{Optimisation}} \\
Optimiser                  & AdamW, $\beta = (0.9, 0.999)$, $\varepsilon = 10^{-8}$ \\
Learning rate              & $10^{-3}$, constant (no warm-up or decay) \\
Weight decay               & $0.01$ \\
Batch size                 & 16 \\
Epochs / gradient steps    & 100 / 1600 \\
Gradient clipping          & Global $\ell_2$ norm $1.0$ \\
Weight EMA, dropout        & None \\
Hardware, wall-clock time  & Single CPU, ${\approx}2$~h \\
\midrule
\multicolumn{2}{l}{\textit{Sampling}} \\
Solver                     & Deterministic Euler, 64 steps \\
Schedule                   & Karras, $\rho = 7$ \\
Initialisation             & Regression baseline $+\ \sigma_{\max}\varepsilon$, $\varepsilon \sim \mathcal{N}(0, I)$ \\
Samples per scenario       & 5 \\
\bottomrule
\end{tabular}
\end{table}